\documentclass[reqno, 11pt]{amsart}
\pdfoutput=1 

\usepackage[usenames, dvipsnames]{xcolor}

\usepackage{amsmath,amssymb,amsthm,mathtools}
\usepackage{tikz}
\usepackage[T1]{fontenc}
\usepackage[utf8]{inputenc}
\usepackage{enumerate}
\usepackage[letterpaper, margin=1in]{geometry}
\usepackage[noBBpl]{mathpazo} 

\usepackage[longnamesfirst,numbers,sort&compress]{natbib}
\makeatletter
\def\NAT@spacechar{~}
\makeatother

\usepackage{thmtools}
\usepackage{thm-restate}
\usepackage{hyperref}
\hypersetup{
colorlinks,
linkcolor={cDeepBlue},
citecolor={green!60!black},
urlcolor={blue!80!black},
pdftitle = {Notes on Aharoni's rainbow cycle conjecture},
pdfauthor = {Tony Huynh and Jackson Goerner and Katie Clinch and Freddie Illingworth}
}
\usepackage[compress, capitalise, nameinlink, noabbrev]{cleveref} 

\declaretheorem[name = Theorem, numberwithin = section, style = plain, refname = {Theorem,Theorems}, Refname = {Theorem,Theorems}]{theorem}
\declaretheorem[name = Claim, numberlike = theorem, style = plain, refname = {Claim,Claims}, Refname = {Claim,Claims}]{claim}

\declaretheorem[name = Lemma, numberlike = theorem, style = plain, refname = {Lemma,Lemmas}, Refname = {Lemma,Lemmas}]{lemma}

\declaretheorem[name = Conjecture, numberlike = theorem, style = plain, refname = {Conjecture,Conjectures}, Refname = {Conjecture, Conjetures}]{conjecture}
\declaretheorem[name = Question, numberlike = theorem, style = plain, refname = {Question,Questions}, Refname = {Question, Questions}]{question}

\newcommand{\rg}{\operatorname{rg}}
\newcommand{\CH}{Caccetta-Häggkvist }

\newcommand{\cO}{\mathcal{O}}

\newcommand{\defn}[1]{\textcolor{cMaroon}{\emph{#1}}}
\newcommand{\df}{\mathsf{def}}
\renewcommand{\geq}{\geqslant}
\renewcommand{\leq}{\leqslant}

\DeclarePairedDelimiter{\abs}{\lvert}{\rvert}
\DeclarePairedDelimiter{\ceil}{\lceil}{\rceil}

\definecolor{cMaroon}{HTML}{93152a}

\definecolor{cYellow}{HTML}{f0e442}
\definecolor{cOrange}{HTML}{f7903b}
\definecolor{cRed}{HTML}{ff3341}
\definecolor{cPurple}{HTML}{da58c2}
\definecolor{cIndigo}{HTML}{8a58da}
\definecolor{cDeepBlue}{HTML}{4a8ce8}
\definecolor{cLightBlue}{HTML}{43e0ef}
\definecolor{cGreen}{HTML}{3ef450}

\tikzset{every node/.style={
    fill=black,
    draw=black,
    text=white,
    inner sep=1pt,
    minimum width=5mm,
    font=\small,
    line width=1.5pt,
    circle
}}
\tikzset{not_vert/.style={
    fill=white,
    text=black,
    draw opacity=0,
    minimum width=0
}}
\tikzset{every path/.style={
    line width=2pt
}}

\begin{document}
\title{Notes on Aharoni's rainbow cycle conjecture}
\author{Katie Clinch}
\author{Jackson Goerner}
\author{Tony Huynh}
\author{Freddie Illingworth}
\address[Katie Clinch]{School of Computer Science and Engineering, UNSW, Sydney, Australia} 
\email{k.clinch@unsw.edu.au} 
\address[Jackson Goerner]{School of Mathematics, Monash University, Melbourne, Australia}
\email{jackson.goerner@monash.edu}
\address[Tony Huynh]{Dipartimento di Informatica, Sapienza Università di Roma, Italy}
\email{huynh@di.uniroma1.it}
\address[Freddie Illingworth]{Mathematical Institute, University of Oxford, Oxford, United Kingdom}
\email{illingworth@maths.ox.ac.uk}
\thanks{Research of F.I.\ supported by EPSRC grant EP/V007327/1.  Research of T.H. was supported by the Australian Research Council.} 

\address[]{}
\email{}
\thanks{}
\begin{abstract}
  In 2017, Ron Aharoni made the following conjecture about rainbow cycles in edge-coloured graphs: If $G$ is an $n$-vertex graph whose edges are coloured with  $n$ colours and each colour class has size at least $r$, then $G$ contains a rainbow cycle of length at most $\ceil{\frac{n}{r}}$.  One motivation for studying Aharoni's conjecture is that it is a strengthening of the Caccetta-H\"{a}ggkvist conjecture on digraphs from 1978.

In this article, we present a survey of Aharoni's conjecture, including many recent partial results and related conjectures.  We also present two new results.  Our main new result is for the $r=3$ case of Aharoni's conjecture.  We prove that if $G$ is an $n$-vertex graph whose edges are coloured with $n$ colours and each colour class has size at least $3$, then $G$ contains a rainbow cycle of length at most $\frac{4n}{9}+7$. We also discuss how our approach might generalise to larger values of $r$.

\end{abstract}
\maketitle

\section{Introduction}
In 1978, Caccetta and H\"{a}ggkvist made the following conjecture about directed cycles in digraphs.  

\begin{conjecture}[Caccetta-Häggkvist] \label{conj:CH}
For all positive integers $n,r$, every simple\footnote{A digraph is \defn{simple} if for all vertices $u$ and $v$, there is at most one arc from $u$ to $v$.} $n$-vertex digraph with minimum outdegree at least $r$ contains a directed cycle of length at most $\ceil{\frac{n}{r}}$.
\end{conjecture}

 Despite considerable effort from numerous researchers, the Caccetta-Häggkvist conjecture remains open. A complete summary of the plethora of results related to the Caccetta-Häggkvist is beyond the scope of this survey. We refer the interested reader to~\citet{sullivan06} for a brief synopsis. We instead focus on the following generalisation of \cref{conj:CH}.
 
 \begin{conjecture}[Aharoni] \label{conj:aharoni}
Let $G$ be a simple\footnote{An edge-coloured graph $G$ is \defn{simple} if each colour class does not contain parallel edges.} edge-coloured graph with $n$ vertices and $n$ colours, where each colour class has size at least $r$.  Then $G$ contains a rainbow\footnote{A subgraph of an edge-coloured graph is \defn{rainbow} if no two of its edges are the same colour.} cycle of length at most $\ceil{\frac{n}{r}}$.
\end{conjecture}

For completeness, we now give a proof that \cref{conj:aharoni} indeed implies \cref{conj:CH}.  In fact, as noted in~\cite{DDFGGHMM21}, the following weakening of Aharoni's conjecture already implies the Caccetta-Häggkvist conjecture.  

\begin{conjecture}[\citet{DDFGGHMM21}] \label{conj:weakaharoni}
Let $G$ be a simple edge-coloured graph with $n$ vertices and $n$ colours, where each colour class has size at least $r$.  Then $G$ contains a cycle $C$ of length at most $\ceil{\frac{n}{r}}$ such that no two incident edges of $C$ are the same colour. 
\end{conjecture}

\begin{proof}[Proof of Conjecture~\ref{conj:CH}, assuming Conjecture~\ref{conj:weakaharoni}]
Let $D$ be a simple digraph with $n$ vertices and minimum outdegree at least $r$. Let $G$ be the graph obtained from $D$ by forgetting the orientations of all arcs.  Colour $uv \in E(G)$ with colour $u$ if $(u,v) \in E(D)$.  Clearly, this colouring uses $|V(D)|=n$ colours.  Moreover, since $D$ has minimum outdegree at least $r$, each colour class has size at least $r$.  Therefore, by Conjecture~\ref{conj:weakaharoni}, $G$ contains a properly edge-coloured cycle $C$ of length at most $\ceil{\frac{n}{r}}$. The set of arcs in $D$ corresponding to the edges of $C$ must be a directed cycle; otherwise $C$ is not properly edge-coloured.  
\end{proof}

Despite the fact that Aharoni's conjecture implies the \CH conjecture, a proof of Aharoni's conjecture may be easier to find than a proof of the \CH conjecture.  Although this might sound counterintuitive, the method of proving a stronger statement is very common in combinatorics. For example,~\citet{thomassen94} found a beautiful short proof that every planar graph is 5-colourable by proving a stronger list colouring version of the theorem.  Moreover, generalisation often leads to new techniques and new questions which one would not even consider in the original setting.  We will see that this is the case for Aharoni's conjecture in the next section.   


\section{Related Results and Conjectures}
In this section, we survey results and conjectures related to rainbow cycles. For a general 
survey on rainbow sets, we refer the reader to~\citet{AB22}.

\subsection{Larger Colour Classes} Much of the research on the Caccetta-Häggkvist conjecture has focused on the directed triangle case ($r=\ceil{\frac{n}{3}}$). A natural strategy is to increase the outdegree condition until one can prove the existence of a directed triangle.  The best result in this direction is the following result of~\citet{HKN17}, which uses the flag algebra machinery developed by~\citet{Razborov13}.

\begin{theorem}
Every simple $n$-vertex digraph with minimum outdegree at least $0.3465n$ contains a directed triangle.  
\end{theorem}

Similarly, for Aharoni's conjecture, one can ask how large must the colour classes be to ensure a rainbow cycle of length at most $\ceil{\frac{n}{r}}$? The first non-trivial bound was proven by~\citet{HPPS21}.

\begin{theorem} \label{thm:rlogr}
Let $r \geq 2$, $n \geq 1$, $G$ be a simple edge-coloured graph with $n$ vertices and $n$ colours, where each colour class has size at least $301r\log r$.  Then $G$ contains a rainbow cycle of length at most $\ceil{\frac{n}{r}}$.  
\end{theorem}

\cref{thm:rlogr} was subsequently improved by~\citet{HS22}, who removed the $\log r$ term.  

\begin{theorem}
Let $r, n \in \mathbb{N}$, $G$ be a simple edge-coloured graph with $n$ vertices and $n$ colours, where each colour class has size at least $10^{11}r$. Then $G$ contains a rainbow cycle of length at most $\frac{n}{r}$.
\end{theorem}

\subsection{Number of Colour Classes} \label{sec:numberofcolours}
One appealing aspect of Aharoni's conjecture is that we can ask what happens when the number of colour classes is different from the number of vertices.  Note that this question does not even make sense in the digraph setting.   To be precise, we now define a function $f(n,t,r)$ which will be useful to state many of the results that appear in this survey. The \defn{rainbow girth} of an edge-coloured graph $G$, denoted $\rg (G)$, is the length of a shortest rainbow cycle in $G$. If $G$ does not contain a rainbow cycle, then $\rg (G)=\infty$. Let 
\[f(n,t,r)\coloneqq \max \{\rg (G)\},\]
where the maximum is taken over all simple edge-coloured graphs $G$ with $n$ vertices and at least $t$ colours, such that each colour class has size at least $r$.  

We can rephrase Aharoni's conjecture via our function $f(n,t,r)$ as follows.

\begin{conjecture}[Aharoni]
For all $n,r \geq 1$,
\[
f(n,n,r) \leq \ceil{\tfrac{n}{r}}.
\]
\end{conjecture}

We believe that there is no reason to restrict attention to the case $t=n$, and that the following question is of independent interest.  

\begin{question} \label{question:f}
Obtain good upper and lower bounds for $f(n,t,r)$ for all $n,t$, and $r$.  
\end{question}

 An important special case of \cref{question:f} is when $r=1$, which was considered by \citet{BS02}.

 \begin{theorem} \label{thm:BS}
For all $n \geq 4$ and $k \geq 2$, 
\[
f(n, n+k, 1) \leq \frac{2(n+k)}{3k} (\log k + \log \log k +4).
\]
 \end{theorem}

 In other words, \cref{thm:BS} asserts that every $n$-vertex graph with at least $n+k$ edges contains a cycle of length at most $\frac{2(n+k)}{3k} (\log k + \log \log k +4)$.  This is a key tool used in many of the results presented in this survey. 

 \citet{HS22} also obtained the following bounds when the number of colours is more than the number of vertices.    

\begin{theorem}
For all $n \geq 1$ and $k \geq 2$, 
\[
f(n, n+k, 10^9k) \leq \min\{6, \frac{n(\log k)^2}{10k^{3/2}} + 14\log k \}.
\]
\end{theorem}

When the number of colours is \emph{less} than the number of vertices, \citet{DDFGGHMM21} obtained the following tight bounds for $r=2$. 
\begin{theorem} \label{thm:sharp}
For all $n \geq 3$ and $t \leq n$, 
\[
f(n, t, 2) = 
\begin{cases}
\infty & \text{if $t \leq n-2$,} \\
n-1 & \text{if $t=n-1$,} \\
\ceil{\frac{n}{2}} &\text{if $t=n$}.
\end{cases}
\]
\end{theorem}

\subsection{Structured Colour Classes}
In the proof of the \CH conjecture (assuming Aharoni's conjecture) given in the Introduction, the colour classes of the derived edge-colouring are all stars.  Therefore, it is natural to ask what happens when the colour classes are not stars.  Note that if a colour class is not a star, then it must contain a matching of size 2, or it is a triangle. 

In the extreme case when all colour classes contain a matching of size 2, \citet{AH21} proved that there is a much shorter rainbow cycle than the $\ceil{\frac{n}{2}}$ bound which follows from~\cite{DDFGGHMM21}.  

\begin{theorem}
There exists an absolute constant $C$ such that if $G$ is a simple edge-coloured graph with $n$ vertices, $n$ colours, and each colour class is a matching of size 2, then $G$ contains a rainbow cycle of length at most $C\log n$.   
\end{theorem}

Kevin Hendrey (private communication) proved that a $\cO(\log n)$
bound also holds in the case when all colour classes are a triangle.  The same proof also appears in a recent paper of~\citet{ABCGZ21}.

\begin{theorem}
There exists an absolute constant $C$ such that if $G$ is a simple edge-coloured graph with $n$ vertices, $n$ colours, and each colour class is a triangle, then $G$ contains a rainbow cycle of length at most $C\log n$.  
\end{theorem}

Some mixed cases were investigated by~\citet{guo22}.

\begin{theorem}
There exists an absolute constant $C$ such that if $G$ is a simple edge-coloured graph with $n$ vertices, $n$ colours, and each colour class is a matching of size 2 or a triangle, then $G$ contains a rainbow cycle of length at most $C\log n$.  
\end{theorem}

\begin{theorem}
For any constants $0 \leq \alpha < 1$ and $0 \leq \beta \leq \alpha$ with $\beta < (1 - \alpha)/3$, there exists a constant $C(\alpha, \beta)$ such that if $G$ is an $n$-vertex simple edge-coloured graph containing at least $(\alpha - \beta)n$ color classes consisting of a single
edge and at least $(1 - \alpha - \beta)n$ color classes consisting of a triangle, then $G$ contains a rainbow cycle of length at most 
$C(\alpha, \beta) \log n$.
\end{theorem}

\subsection{Rainbow Triangles}
The rainbow triangle case ($r = \ceil{n/3}$) of Aharoni's conjecture states that $f(n,n, \ceil{n/3}) \leq 3$. This is of course still open since the directed triangle case of the \CH conjecture is still open.  However, there have been partial results which increase the number of colours or sizes of the colour classes.  Two such results were obtained by~\citet{ADH19}.

\begin{theorem}
For all $n \geq 1$, 
\[
f(n, 9n/8, n/3) \leq 3.
\]
\end{theorem}

\begin{theorem}For all $n \geq 1$, 
\[
f(n, n, 2n/5) \leq 3.
\]
\end{theorem}

Both of these results have been subsequently improved by~\citet{HQS22}.

\begin{theorem}
For all $n \geq 1$, 
\[
f(n, 1.1077n, n/3) \leq 3.
\]
\end{theorem}

\begin{theorem}
For all $n \geq 1$,
\[
f(n, n, 0.3988n) \leq 3.
\]
\end{theorem}

\citet{ADGMS20} proved the following theorem showing when a simple edge-coloured graph with 3 colours contains a rainbow triangle. Their theorem actually implies Mantel's theorem.  

\begin{theorem} \label{thm:rainbowmantel}
Every $n$-vertex simple edge-coloured graph with 3 colours and each colour class of size at least  $\frac{26-2\sqrt{7}}{81} \cdot n^2$ contains a rainbow triangle.
\end{theorem}

They also prove that the constant $\frac{26-2\sqrt{7}}{81}$ in~\cref{thm:rainbowmantel} is best possible. However, their extremal example has rainbow 2-cycles.  Our first new result is the following sharp thresholds for $f(n,3,r)$.

\begin{theorem} \label{thm:sharp}
For all $n \geq 100$, 
\[
f(n, 3, r) = 
\begin{cases}
\infty & \text{if $r \leq \lfloor \binom{n}{2} /3 \rfloor$,} \\
2 & \text{if $\lfloor \binom{n}{2} /3 \rfloor < r \leq \binom{n}{2}$,} \\
1 & \text{if $r > \binom{n}{2}$.}
\end{cases}
\]
\end{theorem}

\begin{proof}
First suppose $r > \lfloor \binom{n}{2} /3 \rfloor$. Let $G$ be a simple edge-coloured graph with $n$ vertices, 3 colours and each colour class of size at least $r$.  Since $K_n$ has only $\binom{n}{2}$ edges, $G$ must contain a loop or a rainbow 2-cycle. Moreover, if  $r > \binom{n}{2}$, then $G$ must contain a loop.  The $n$-vertex edge-coloured graph with a red, blue, and green edge between every pair of vertices proves equality.  

We now show that $f(n, 3, r) = \infty$ if $n \geq 100$ and $r \leq \lfloor \binom{n}{2} /3 \rfloor$.  Partition the vertices of $K_n$ as $X \cup Y \cup Z$ where $|X|=|Y|=\lceil \frac{2n}{5} \rceil$.  Colour all edges in $K_n[X]$ or $ K_n[Y]$ or $K_n[Z]$ red, all edges between $X$ and $Y$ blue, and all edges between $Z$ and $X \cup Y$ green.  Observe that this is an edge-colouring of $K_n$ with no rainbow triangle.  The colour classes do not quite have the same size, which we now fix.  
Since $n \geq 100$, we have $|X||Y| \leq \binom{n}{2} /3$ and $|Z|(|X|+|Y|) \leq \binom{n}{2} /3$.  Thus, we may recolour some of the edges in $K_n[X]$ blue, and some of the edges in $K_n[Y]$ green, so that the number of red, blue, and green edges differ by at most 1. This new edge-colouring of $K_n$ still does not contain a rainbow triangle. This of course implies that there are no rainbow cycles since there are only 3 colours.  
\end{proof}

\subsection{Small Values of \texorpdfstring{$r$}{r}}
The \CH conjecture is known to hold for small values of $r$.  The $r=1$ case is trivial, the $r=2$ case was proven by~\citet{CH78}; the $r=3$ case was proven by~\citet{hamidoune87}; and the $r \in \{4,5\}$ cases were proven  by~\citet{HR87}.  For Aharoni's conjecture, the $r=2$ case was proven by~\citet{DDFGGHMM21}.

\begin{theorem} \label{thm:r=2}
For all $n \geq 1$,
\[
f(n,n,2) \leq \ceil{\tfrac{n}{2}}.
\]
\end{theorem}

In~\cref{sec:main}, we will prove the following result for $r=3$.

\begin{theorem} \label{thm:r=3}
For all $n \geq 1$,
\[
f(n,n,3) \leq \tfrac{4n}{9} +7.
\]
\end{theorem}



\subsection{Non-uniform Versions}
For the \CH conjecture, it is natural to ask if there is a version which takes into account \emph{all} the outdegrees rather than just the minimum outdegree. Seymour (see~\cite{hompe19}) proposed the following generalization. Given a digraph $D$ with no sink\footnote{A \defn{sink} is a vertex with outdegree zero.}, define
\[
\psi(D):=\sum_{v \in V(D)} \frac{1}{\deg^+(v)}.
\]

\begin{conjecture} \label{seymour}
Every simple digraph $D$ with no sink contains a directed cycle of length at most $\lceil \psi(D) \rceil$.
\end{conjecture}

Note that in the case that all outdegrees are $r$, then $\psi(D)=\frac{n}{r}$, so~\cref{seymour} implies the \CH conjecture.  Unfortunately, \cref{seymour} was disproved by~\citet{hompe19}.  However,~\citet{ABCGZ21} proved that "half" of~\cref{seymour} holds. 

\begin{theorem} \label{almostseymour}
   Every simple digraph $D$ with no sink contains a directed cycle of length at most $2 \psi(D) $.
 \end{theorem}

\cref{almostseymour} has a natural generalization in the rainbow setting.  Given an edge-coloured graph $G$, define
\[
\psi(G)=\sum_{A} \frac{1}{|A|},
\]
where the sum is taken over all colour classes $A$ of $G$. 

\begin{conjecture} \label{nonuniformrainbow}
   Every simple edge-coloured graph $G$ with the same number of colours as vertices contains a rainbow cycle of length at most $2\psi(G)$.
\end{conjecture}

In the case that all colour classes have size at most 2,~\citet{ABCGZ21} proved the following strengthening of~\cref{nonuniformrainbow}. 

\begin{theorem} \label{2strong}
   Let $G$ be a simple edge-coloured graph with the same number of colours as vertices and such that each colour class has size at most 2. Then $G$ contains a rainbow cycle of length at most $\lceil \psi(G) \rceil$.
\end{theorem}

Note that~\cref{2strong} is a strengthening of~\cref{thm:r=2} since it allows colour classes of size 1.

\subsection{Matroids} We can generalise Aharoni's conjecture to any setting in which the notion of `cycle' makes sense.  One natural candidate is that of a \emph{matroid}.   For the reader unfamiliar with matroids, we introduce all the necessary definitions now. For a more thorough introduction to matroids, we refer the reader to Oxley~\cite{oxley11}. 

A \defn{matroid} is a pair $M=(E, \mathcal C)$ where $E$ is a finite set, called the \defn{ground set} of $M$, and $\mathcal C$ is a collection of subsets of $E$, called \defn{circuits}, satisfying
\begin{enumerate}
    \item $\emptyset \notin \mathcal C$, 
    \item if $C'$ is a proper subset of $C \in \mathcal C$, then $C' \notin \mathcal C$,
    \item if $C_1$ and $C_2$ are distinct members of $\mathcal C$ and $e \in C_1 \cap C_2$, then there exists $C_3 \subseteq (C_1 \cup C_2) \setminus \{e\}$.
\end{enumerate}

We now give some examples of matroids.  Let $G$ be a graph.  We will consider two different matroids with ground set $E(G)$.  The circuits of the first matroid are the (edges of) cycles of $G$. This is the \defn{cycle matroid} of $G$, denoted $M(G)$.   A matroid is \defn{graphic} if it is isomorphic to the cycle matroid of some graph. A \defn{cocycle} of $G$ is an inclusion-wise minimal edge-cut of $G$.  The collection of cocycles of $G$ is also a matroid, called the \defn{cocycle matroid} of $G$, and is denoted $M(G)^*$.  A matroid is \defn{cographic} if it is isomorphic to the cocycle matroid of some graph. 

Let $\mathbb F$ be a field.  An \defn{$\mathbb F$-matrix} is a matrix with entries in $\mathbb F$.  Let $A$ be an $\mathbb F$-matrix whose columns are labelled by a finite set $E$. The \defn{column matroid} of $A$, denoted $M[A]$, is the matroid with ground set $E$ whose circuits correspond to the minimal (under inclusion) linearly dependent columns of $A$.  A matroid is \defn{representable over $\mathbb{F}$} if it is isomorphic to $M[A]$ for some $\mathbb{F}$-matrix $A$.  A matroid is \defn{binary} if it is representable over the two-element field,  and it is \defn{regular} if it is representable over every field.  

In order to formulate Aharoni's conjecture for matroids, we need to define simple matroids and how to express the number of vertices of a graph as a matroid parameter.  We do this now.  A matroid is \defn{simple} it it does not contain any circuits of size $1$ or $2$.  A set $I \subseteq E$ is \defn{independent} if it does not contain a circuit.  
The \defn{rank} of $X \subseteq E$ is the size of a largest independent set contained in $X$, and is denoted $r_M(X)$.  The \defn{rank} of $M$ is $r(M)\coloneqq r_M(E)$.  Notice that the number of vertices of a connected graph $G$ is $r(M(G))-1$. Thus, Aharoni's conjecture can be phrased in the language of matroids as follows.

\begin{conjecture}[Aharoni] \label{conj:matroid}
Let $M$ be a simple rank-$(n-1)$ graphic matroid and $c$ be a colouring of $E(M)$ with $n$ colours, where each colour class has size at least $r$. Then $M$ contains a rainbow circuit of size at most $\ceil{\frac{n}{r}}$.
\end{conjecture}

One way to generalise \cref{conj:CH} is to replace `graphic matroid' with some larger superclass of matroids.  In~\cite{DDFGGHMM21}, it was shown that one cannot replace `graphic matroid' by `binary matroid' in \cref{conj:matroid}.
\begin{theorem}
For all  $n \geq 6$, there exists a simple rank-$(n-1)$ binary matroid $M$ on $2n$ elements, and a colouring of $E(M)$ where each colour class has size $2$, such that all rainbow circuits of $M$ have size strictly greater than $\lceil \frac{n}{2} \rceil$.
\end{theorem}

The main result of~\cite{DDFGGHMM21} can be phrased in matroid language as follows.  

\begin{theorem} \label{thm:graphicmatroid}
Let $M$ be a simple rank-$(n-1)$ graphic matroid and $c$ be a colouring of $E(M)$ with $n$ colours, where each colour class has size at least $2$. Then $M$ contains a rainbow circuit of size at most $\ceil{\frac{n}{2}}$.
\end{theorem}

In~\cite{DDFGGHMM21}, it is also proved that the matroid analogue of \cref{thm:r=2} holds for cographic matroids.

\begin{theorem} \label{thm:cographic}
Let $N$ be a simple rank-$(n-1)$ cographic matroid and $c$ be a colouring of $E(N)$ with $n$ colours, where each colour class has size at least $2$. Then $N$ contains a rainbow circuit of size at most $\ceil{\frac{n}{2}}$.
\end{theorem}

Regular matroids are a well-studied superclass of graphic matroids, and are `essentially' graphic or cographic via Seymour's regular matroid decomposition theorem~\cite{Seymour80}.  Therefore, by combining~\cref{thm:graphicmatroid},~\cref{thm:cographic}, and Seymour's regular matroid decomposition theorem, it may be possible to prove the following conjecture.    

\begin{conjecture}[\citet{DDFGGHMM21}] \label{conj:regular}
Let $M$ be a simple rank-$(n-1)$ regular matroid and $c$ be a colouring of $E(M)$ with $n$ colours, where each colour class has size at least $2$. Then $M$ contains a rainbow circuit of size at most $\ceil{\frac{n}{2}}$.
\end{conjecture}

Of course, we can also consider~\cref{question:f} for various classes of matroids.  \citet{BS22} obtained one such result. 

\begin{theorem}
Let $M$ be an $n$-element rank-$t$ binary matroid whose ground set is coloured with $t$ colours.  Then $M$ either contains a rainbow circuit or a monochromatic cocircuit.\footnote{A \defn{monochromatic cocircuit} is a circuit in the dual matroid whose elements are all the same colour.}  
\end{theorem}

\citet[Theorem 4]{BS22} also show that it is possible to characterize binary matroids as exactly those matroids which do not admit a specific type of colouring with no rainbow circuits.  They also prove that if a simple graph $G$ has an edge-colouring with no rainbow cycle, where each colour class has size at most 2, then $G$ is independent in the 2-dimensional rigidity matroid.\footnote{Equivalently, for every $X \subseteq V(G)$ with $|X| \geq 2$, $|E(G[X])| \leq 2|X|-3$. }

\section{Proof of Main Theorem} \label{sec:main}
In this section, we prove our main theorem.

\begin{theorem} \label{thm:main}
Every simple edge-coloured graph with $n$ vertices, $n$ colours, and each colour class of size at least 3, contains a rainbow cycle of length at most $\tfrac{4n}{9}+7$.  
\end{theorem}

\subsection{Excess-\texorpdfstring{$k$}{k} Graphs}
We begin by establishing some basic properties about graphs which have at least $k$ more edges than vertices, which we call \defn{excess-$k$ graphs}.  The first property follows from ~\cref{thm:BS}: for all $n \geq 4$ and $k \geq 2$, every $n$-vertex, excess-$k$ graph has a cycle of length at most $\frac{2(n+k)}{3k} (\log k + \log \log k +4)$.  We require the following tighter bounds when $k \leq 2$.  

\begin{lemma} \label{claim:excess1}
Every $n$-vertex, excess-1 graph has a cycle of length at most $\frac{2n}{3}+1$.
\end{lemma}

The proof is easy and is omitted (see~\cite{DDFGGHMM21} for a proof of a stronger claim).  

\begin{lemma} \label{claim:excess2}
Every $n$-vertex, excess-2 graph $H$ has a cycle of length at most $\frac{n}{2}+1$.
\end{lemma}

\begin{proof}
 A \defn{block} of $H$ is a 2-connected subgraph of $H$ which is maximal under the subgraph relation.  If $H$ contains at least two blocks, then $H$ contains two cycles which meet in at most one vertex.  Hence, one of these cycles has length at most $\frac{n}{2}+1$.  So, we may assume $H$ contains exactly one block $B$.  Note that $B$ is excess-2, since $H$ is excess-2.  Let $C \cup P_1 \cup \dots \cup P_k$ be an ear-decomposition of $B$.  Note that $k \geq 2$, since $B$ is excess-2.  Let $B'=C \cup P_1 \cup P_2$. 
 Either $B'$ contains two cycles which meet in at most two vertices, or $B'$ is a subdivision of $K_4$.  In the first case, one of the two cycles has length at most $\frac{n}{2}+1$.  In the second case, $B'$ contains four cycles $C_1, \dots, C_4$ such that $\sum_{i\in [4]}|V(C_i)| = 2|V(B')|+4$. Thus, one of these four cycles has length at most $\frac{|V(B')|}{2}+1 \leq \frac{n}{2}+1$.  
 \end{proof}

 We now prove that the largest stable set of an excess-$k$ graph is roughly half the number of vertices.

\begin{lemma} \label{maxstableset}
Let $H$ be a simple excess-$k$ graph of minimum degree at least 2. Then a maximum stable set of $H$ has size at most $\frac{|V(H)|+k}{2}$.  
\end{lemma}

\begin{proof}
We proceed by induction on $|V(H)|+|E(H)|$.  We may assume that $H$ has exactly $k$ more edges than vertices.  If $H$ is the disjoint union of $H_1, \dots, H_\ell$ and $H_i$ has excess $k_i$, then $k_1 + \dots + k_\ell=k$.  Therefore, we are done by applying induction to each $H_i$. 
 So, we may assume $H$ is connected.  If $k=0$, then $H$ is a cycle, so the lemma clearly holds.  Now suppose $k \geq 1$.  Let $X$ denote the set of vertices of $H$ of degree at least 3. If $|X|=1$, then $H$ consists of $k+1$ cycles which meet at the same vertex.  It is easy to see that the lemma holds in this case.  So, we may assume $|X| \geq 2$. Suppose $u,v \in X$ and $uv \in E(G)$.  By induction, a maximum stable set of $H - uv$ has size at most $\frac{|V(H)|+k-1}{2}$.  Hence, a maximum stable set of $H$ has size at most $\frac{|V(H)|+k-1}{2}$.  
So, we may assume that no two vertices of $X$ are adjacent.  

Let $P$ be the shortest path between any two vertices of $X$.  By the minimality of $P$, each internal vertex of $P$ has degree 2.  Let $H'=H-I$, where $I$ is the set of internal vertices of $P$.  Note that $H$' has excess $k-1$ and minimum degree at least 2.  Let $S$ be a maximum size stable set in $H$.  Note that $S \cap V(H')$ is a stable set in $H'$.  By induction, $|S \cap V(H')| \leq (|V(H')|+k-1)/2$.  Also, $|S \cap I| \leq (|I|+1)/2$.  Thus, $|S| = |S \cap V(H')|+|S \cap I| \leq (|V(H)|+k)/2$.
\end{proof}

We finish by establishing the following lemma about `minimal' excess-$k$ rainbow subgraphs of an edge-coloured graph.



 \begin{lemma} \label{fewchords}
Let $H$ be a simple edge-coloured graph where each colour class has size at most $r$.  Let $R$ be an excess-$k$ rainbow subgraph of $H$ such that $V(R)$ is minimal under inclusion. Then $H$ contains a rainbow 2-cycle, or there are at most $\max\{\binom{2k+2}{2},6k(r-1)\}$ chords of $R$ in $H$.  
\end{lemma}

\begin{proof}
We may assume that $H$ does not contain parallel edges; otherwise $H$ contains a rainbow 2-cycle.  By the minimality of $R$, $\deg_R(v) \geq 2$ for all $v \in V(R)$. Let $V(R)=X \cup Y$, where $X$ is the set of vertices of degree at least 3, and $Y$ is the set of degree-2 vertices of $R$. By the Handshaking Lemma, $|X| \leq 2k$.  

 Say that a chord of $R$ is \defn{novel} if its colour does not appear in $R$, and \defn{plain} otherwise. Suppose $e=uv$ is a novel chord.  If $|Y| \geq 3$, then there is a vertex $y \in Y$ such $y \notin \{u,v\}$, and so $(R \cup e) - y $  contradicts the minimality of $R$.  Thus, there are either no novel chords, or the total number of chords is at most $\binom{|X|+|Y|}{2} \leq \binom{2k+2}{2}$. 
 
 We may thus assume that there are no novel chords.  Let $e \in E(H)$ be a plain chord of $R$.  Since $e$ is plain, there is an edge $f \in E(R)$ of the same colour as colour $e$.  If both ends of $f$ are in $Y$, then $R'\coloneqq (R \cup e) \setminus f$ contains a degree-1 vertex $x$.  But now $R'-x$ contradicts the minimality of $R$.  Therefore, at least one end of $f$ is in $X$.  Let $R_0$ be the multigraph on vertex set $X$ obtained by suppressing all degree-2 vertices.  Note that $R_0$ also has excess $k$ and at most $2k$ vertices. Therefore, $|E(R_0)| \leq 3k$.  There are at most $2|E(R_0)|$ edges of $R$ which are incident to a vertex in $X$.  Each of these edges can correspond to at most $r-1$ plain chords.  Hence, there are at most $6k(r-1)$ plain chords.
\end{proof}

\subsection{The set-up} Let $G$ be a simple edge-coloured graph with $n$ vertices, $n$ colours, and each colour class of size exactly 3.  An \defn{$r$-star} is a star with $r$ edges.  A colour class of $G$ is a \defn{star class} if it is a 3-star. A vertex of $G$  is  a \defn{star vertex} if it is the centre of a star class, and is otherwise a \defn{non-star vertex}. Let $S$ denote the set of star vertices of $G$, and $N$ denote the set of non-star vertices of $G$.  Let $N'$ be the set of non-star classes.  Since star classes may be centred at the same vertex, we have $|N'| \leq |N|$.  We will do a case analysis depending on whether $|N| \geq 8$ or $|N| \leq 7$.  

\subsection{Many non-star vertices}  Throughout this section we suppose $|N| \geq 8$.  Let $\{x,y\} \subseteq N$.  We say that a colour class $A$ \defn{dominates} $\{x,y\}$ if every edge in $A$ has an end in $\{x,y\}$.  

\begin{claim} \label{claim:nodominate}
There exists $\{x,y\} \subseteq N$ such that no colour class dominates $\{x,y\}$.  
\end{claim}

\begin{proof}
Suppose $A$ is a star class.  Then $A$ is a 3-star centred at a vertex not in $N$.  Therefore, for all $\{u,v\} \subseteq N$, $A$ does not dominate $\{u,v\}$ since at least one leaf of $A$ is not in $\{u,v\}$.  

Suppose $A$ is a non-star class. Up to isomorphism, $A \in \{K_3, P_3, P_2 \sqcup P_1, P_1 \sqcup P_1 \sqcup P_1\}$, where $P_i$ is a path with $i$ edges, and $\sqcup$ denotes disjoint union.  Let $\gamma(A)$ be the number of vertex covers of $A$ of size 2.  
Observe that the number of pairs $\{u,v\} \subseteq N$ which $A$ dominates is exactly equal to $\gamma(A)$.  We have $\gamma(K_3)=3, \gamma(P_3)=1, \gamma(P_2 \sqcup P_1)=2$, and $\gamma(P_1 \sqcup P_1 \sqcup P_1)=0$.  Thus, every non-star class dominates at most 3 pairs of vertices in $N$.   So, the number of pairs dominated by non-star classes is at most $3|N'| \leq 3|N|< \binom{N}{2}$, since $|N| \geq 8$. Therefore, at least one pair  $\{x,y\} \subseteq N$ is not dominated by any colour class.  
\end{proof}

\begin{claim}
$G$ contains an excess-2 rainbow subgraph $R$.
\end{claim}

\begin{proof}
A \defn{transversal} of $G$ is a subgraph consisting of exactly one edge of each colour.  By~\cref{claim:nodominate}, $G$ has a transversal $R$ such that at least two vertices $x,y$ are not in $V(R)$.  Since $R$ has exactly $n$ edges, $R$ is excess-2 (and clearly rainbow).   
\end{proof}

We now choose $R$ to be an excess-2 rainbow subgraph of $G$ such that $V(R)$ is minimal under inclusion.  

\begin{claim}
$G$ contains a rainbow cycle $C'$ such that $E(R) \cap E(C')=\emptyset$.  
\end{claim}

\begin{proof}
Choose a transversal $R'$ which is edge-disjoint from $R$.  Since $R'$ has $n$ vertices and $n$ edges, $R'$ contains a cycle $C'$.  Clearly, $C'$ is rainbow since $R'$ is rainbow.  
\end{proof}

\begin{claim}
$G$ contains a rainbow cycle $C$ of length at most $\frac{2n}{5} + 7$. 
\end{claim}

\begin{proof}
Let $n_1=|V(R) \setminus V(C')|$, $n_2=|V(R) \cap V(C')|$, and $n_3=|V(C') \setminus V(R)|$.  First suppose $n_3 \geq \frac{n}{5}-12$. Then, $|V(R)| \leq \frac{4n}{5}+12$.  By~\cref{claim:excess2}, $R$ contains a rainbow cycle of length at most $\frac{1}{2} \cdot(\frac{4n}{5}+12)+1 = \frac{2n}{5}+7$.  

Thus, we may assume that $n_3 < \frac{n}{5}-12$. Let $A$ be the subset of edges of $C'$ which are chords of $R$.  Applying~\cref{fewchords} to $R$ in $C' \cup R$ (so $r=2$), we have $|A| \leq 15$. 
Note that $V(R) \cap V(C')$ is a stable set of $C' \setminus A$.  The maximum stable set of $C'$ has size at most $\frac{|V(C')|}{2}$. Deleting one edge from a graph can increase the size of a maximum stable set by at most 1.  Therefore, we conclude that $n_2 \leq n_3+2|A| \leq n_3+30$. 
Thus, 
\[
|V(C')| = n_2+n_3 \leq 2n_3+30 < 2(\tfrac{n}{5}-12)+30 = \tfrac{2n}{5}+6, 
\]
and so we may take $C=C'$ in this case.    
\end{proof}

This completes the case $\abs{N} \geq 8$, since we have found a rainbow cycle of length at most $\frac{2n}{5} + 7$, which is better than the bound of $\frac{4n}{9}+7$ required by~\cref{thm:main}.  

\subsection{Few non-star vertices} We complete the proof by considering the case $\abs{N} \leq 7$.  

\begin{claim} \label{claim:nonstar}
At least one vertex of $G$ is a non-star vertex.  
\end{claim}

\begin{proof}
Suppose every vertex of $G$ is a star vertex.  Since $G$ has the same number of vertices as colours, this implies that at each $v \in V(G)$ there is a star class $S_v$ centred at $v$. Let $D$ be the digraph obtained from $G$ by orienting the edges of $S_v$ away from $v$ for all $v \in V(G)$.  By the $r=3$ case of the \CH conjecture~\cite{hamidoune87}, $D$ contains a directed cycle $\vec{C}$ of length at most $\ceil{\frac{n}{3}}$. Note that $\vec{C}$ corresponds to a rainbow cycle $C$ in $G$.  
\end{proof}

Fix a non-star vertex $z \in V(G)$.  Since $z$ is a non-star vertex, for every colour $a$, there is an edge $e_a$ coloured $a$ such that $e_a$ is not incident to $z$.  Thus, there is a transversal $R$ of $G$ such that $z \notin V(R)$.  In particular, $R$ is an excess-1 rainbow subgraph of $G$ such that $z \notin V(R)$.  Let $R_1$ be an excess-1 rainbow subgraph of $G -z$ such that $V(R_1)$ is minimal under inclusion.  

\begin{claim} \label{claim:7edges}
There exists a transversal $R_2'$ of $G$ such that $z \notin V(R_2')$ and $|E(R_1) \cap E(R_2')| \leq 7$.  
\end{claim}

\begin{proof}
Let $A$ be a star colour class.  Observe that there are at least 2 edges of $A$ which are not incident to $z$.  Thus, there exists an edge $e_A \in A$ such that $e_A \notin E(R_1)$ and $e_A$ is not incident to $z$.  If $A$ is a non-star colour class, there is an edge $e_A \in A$ which is not incident to $z$.  Let $R_2'=\bigcup_A \{e_A\}$, where the union is over all colour classes.  Since there are at most $|N'| \leq |N| \leq 7$ non-star colour classes, $|E(R_1) \cap E(R_2')| \leq 7$.
\end{proof}

Since $z \notin  V(R_2')$, $R_2'$ is excess-1 (and rainbow).  Let $R_2 \subseteq R_2'$ be excess-1 and rainbow such that $V(R_2)$ is minimal under inclusion.  In particular, $R_2$ has minimum degree at least 2.  

\begin{claim} \label{claim:final}
$G$ contains a rainbow cycle $C$ of length at most $\frac{4n}{9}+7$. 
\end{claim}

\begin{proof}
Let $n_1=|V(R_1) \setminus V(R_2)|$, $n_2=|V(R_1) \cap V(R_2)|$, and $n_3=|V(R_2) \setminus V(R_1)|$.  First suppose $n_3 \geq \frac{n}{3}-9$.  Thus, $|V(R_1)| \leq \frac{2n}{3}+9$.  By~\cref{claim:excess1}, $R_1$ contains a rainbow cycle of length at most $\frac{2}{3} \cdot (\frac{2n}{3}+9)+1 = \frac{4n}{9}+7$.  

Thus, we may assume that $n_3 < \frac{n}{3}-9$. Let $A$ be the subset of edges of $R_2$ which are chords or edges of $R_1$.  By~\cref{fewchords} applied to $R_1$ in $R_1 \cup R_2$ (so $r=2$), and~\cref{claim:7edges}, $|A| \leq 6+7=13$.  Note that $V(R_1) \cap V(R_2)$ is a stable set of $R_2 \setminus A$.  By~\cref{maxstableset}, the maximum stable set of $R_2$ has size at most $\frac{|V(R_2)|+1}{2}$. Deleting one edge from a graph can increase the size of a maximum stable set by at most 1.  Therefore, we conclude that $n_2 \leq n_3+2|A|+1 \leq n_3+27$.  By~\cref{claim:excess1}, $R_2$ contains a rainbow cycle of length at most 
\[
\tfrac{2}{3} \cdot (n_2+n_3)+1 \leq \tfrac{2}{3} \cdot (2n_3+27)+1 <  \tfrac{2}{3} \cdot (2(\tfrac{n}{3}-9)+27)+1 = \tfrac{4n}{9}+7. \qedhere
\]
\end{proof}

\cref{claim:final} completes the proof of~\cref{thm:main}.

\section{Generalising our Approach}
 For the $r=3$ case of Aharoni's conjecture, we proved that there is a rainbow cycle of length at most $\frac{4n}{9}+7$.  The additive constant can be easily improved since there are stronger versions of~\cref{fewchords} for $k \leq 2$.  However, we opted to prove~\cref{fewchords} for all $r$ and $k$ in case anyone would like to generalize our approach.
 
 Our strategy was to handle the cases of many non-star vertices and few non-star vertices separately.  For $r=3$, `many' meant at least 8, but for larger $r$ it can be a constant depending on $r$.  Recall, that in the many non-star case, we actually proved a better bound of $\frac{2n}{5}+7$.  Thus, if we can improve the few non-star case, we would get a better theorem. One tempting strategy for the few non-star case is to use the following  `defect' version of the \CH conjecture due to~\citet{shen2000}.

 Let $D$ be a digraph.  For each $r \in \mathbb{N}$, we define the \defn{defect} of $D$ to be
\[
\df_r(D)\coloneqq \sum_{u \in U} (r-\deg^+(u)),
\]
where $U$ is the set of vertices of $D$ of outdegree at most $r$.  

\begin{theorem}[\cite{shen2000}]\label{thm:shen}
Let $D$ be a simple $n$-vertex digraph with no sink, and let $g$ be the length of a shortest directed cycle of $D$.  If $g \geq 2r-1$, then $n \geq r(g-1)+1-\df_r(D)$.  
\end{theorem}

The idea is that if the edges of each star colour class are oriented away from their centres, then we obtain a digraph with small defect.  Unfortunately, the non-star vertices are sinks of this digraph, so we cannot apply~\cref{thm:shen} directly. 
 It may be possible to apply~\cref{thm:shen} in some auxiliary digraph, but we could not see a way to do this. For example, the case that there is exactly one non-star vertex corresponds to the following conjecture about digraphs, which we do not know how to solve. 

 \begin{conjecture}
     Let $D$ be a simple $n$-vertex digraph where each vertex except one has outdegree at least $r$.  Let $\{s_1, t_1\}, \dots, \{s_r, t_r\}$ be $r$ distinct pairs of vertices in $D$.  Then $D$ either contains a directed cycle of length at most $\ceil{\frac{n}{r}}$ or there exist $i \in [r]$, $z \in V(D)$, a directed path $P$ from $s_i$ to $z$, and a directed path $Q$ from $t_i$ to $z$ 
     such that $|E(P)|+|E(Q)|+1 \leq \ceil{\frac{n}{r}}$.
 \end{conjecture}

 For the many non-star case, our method can produce almost edge-disjoint, rainbow subgraphs $H_1, \dots, H_\ell$, each of excess-$k$, where $\ell$ and $k$ depend on the definition of `many' and on $r$.  We conjecture that $H_1 \cup \dots \cup H_\ell$ should contain a short rainbow cycle, but were unable to find a proof.  Our attempted strategy was to introduce a variable for every `region' of the Venn diagram corresponding to $V(H_1), \dots, V(H_\ell)$, and to write out a linear program with constraints coming from~\cref{thm:BS}, \cref{maxstableset}, and~\cref{fewchords}. Unfortunately, it turns out that our LP has a large dual solution and so it cannot prove non-trivial bounds. However, it may be possible to add additional constraints to our LP to obtain good bounds for general $r$.

\subsection*{Acknowledgements} This research was initiated at the \href{https://www.matrix-inst.org.au/events/structural-graph-theory-downunder-ll/}{Structural Graph Theory Downunder II} workshop at the Mathematical Research Institute MATRIX (March 2022). We thank Ron Aharoni and Zach Hunter for carefully reading an earlier version of this paper and making some very helpful suggestions.

  \let\oldthebibliography=\thebibliography
  \let\endoldthebibliography=\endthebibliography
  \renewenvironment{thebibliography}[1]{%
    \begin{oldthebibliography}{#1}%
      \setlength{\parskip}{0ex}%
      \setlength{\itemsep}{0ex}%
  }{\end{oldthebibliography}}

{\fontsize{11pt}{12pt}
\selectfont
\bibliographystyle{DavidNatbibStyle}
\bibliography{references}
}
\printindex
\end{document}